\documentclass[12pt]{article}
\usepackage{latexsym,amsfonts,amssymb,mathrsfs}
\usepackage[dvips]{color}
\usepackage{colordvi,multicol}

\def \cal{\mathcal}

\newtheorem{thm}{Theorem}[section]
\newtheorem{cor}[thm]{Corollary}
\newtheorem{lem}[thm]{Lemma}
\newtheorem{pro}[thm]{Proposition}
\newtheorem{defi}[thm]{Definition}
\newtheorem{rem}[thm]{Remark}

\begin{document}
\title{\bf  Convergences of  Random Variables under Sublinear Expectations}
\author{Ze-Chun Hu\thanks{Corresponding author.\ E-mail addresses: zchu@scu.edu.cn (Z.-C. Hu), 15266479708@163.com (Q.-Q. Zhou).}\\
 {\small College of Mathematics, Sichuan  University, Chengdu, China}\\ \\
 Qian-Qian Zhou\\
 {\small Department of
Mathematics, Nanjing University, Nanjing,  China} }
\date{}
\maketitle



\vskip 0.5cm \noindent{\bf Abstract}\quad In this note, we will survey the existing  convergence results for random variables under sublinear expectations, and prove some new results. Concretely, under the assumption that the sublinear expectation has the monotone continuity property,  we will prove that $L^p$ convergence is stronger than convergence in capacity, convergence in capacity is stronger than convergence in distribution, and give some equivalent characterizations of convergence in distribution. In addition, we give a
dominated convergence  theorem under sublinear expectations, which may have its own interest.

\smallskip


\smallskip


\tableofcontents

\section{Introduction}

It is well known that limit theory plays an important role in probability theory and statistics. Let $(\Omega,\mathcal{F},P)$ be a probability space and $\{X,X_n,n\geq 1\}$ be a sequence of random varialbes. Then we have the following convergences:

\begin{itemize}

\item $\{X_n,n\geq 1\}$  is said to almost surely converge to $X$, if there exists a set $N\in \mathcal{F}$ such that $P(N)=0$ and  $\forall\omega\in \Omega\backslash N, \lim_{n\to\infty}X_n(\omega)=X(\omega)$, which is denoted by $X_n\stackrel{a.s.}{\longrightarrow} X$ or $X_n\to X\ a.s.$

\item $\{X_n,n\geq 1\}$  is said to converge to $X$ in probability, if for any $\varepsilon>0$, $\lim_{n\to\infty}P(\{|X_n-X|\geq \varepsilon\})=0$,  which is denoted by $X_n\stackrel{P}{\longrightarrow} X$.

\item $\{X_n,n\geq 1\}$  is said to $L^p$ converge to $X$ $(p>0)$ if $\lim_{n\to\infty}E[|X_n-X|^p]=0$, which is denoted by $X_n\stackrel{L^p}{\longrightarrow} X$.

\item $\{X_n,n\geq 1\}$  is said to converge to $X$ in distribution, if for any bounded continuous function $f$, $\lim_{n\to\infty}E[f(X_n)]=E[f(X)]$,  which is denoted by $X_n\stackrel{d}{\longrightarrow} X$.

\item $\{X_n,n\geq 1\}$  is said to completely converge to $X$, if for any $\varepsilon>0$, $\sum_{n=1}^{\infty}P(\{|X_n-X|\geq \varepsilon\})<\infty$,  which is denoted by $X_n\stackrel{c.c.}{\longrightarrow} X$ (see \cite{HR47}).

\item $\{X_n,n\geq 1\}$  is said to s-$L^r$ converge to $X$ $(r>0)$ if $\sum_{n=1}^{\infty}E[|X_n-X|^r]<\infty$, which is denoted by $X_n\stackrel{s\mbox{-}L^r}{\longrightarrow} X$ (see \cite[Definition 1.4]{LH15}).
\end{itemize}
Then we have
\begin{eqnarray*}
\begin{array}{ccccccccc}
X_n\stackrel{s\mbox{-}L^r}{\longrightarrow} X &\Rightarrow & X_n\stackrel{c.c.}{\longrightarrow} X& \Rightarrow & X_n\stackrel{a.s.}{\longrightarrow} X &\Rightarrow & X_n\stackrel{P}{\longrightarrow} X &\Rightarrow & X_n\stackrel{d}{\longrightarrow} X,\\
&&&&&&\Uparrow&&\\
&&&&&&X_n\stackrel{L^p}{\longrightarrow} X&&
\end{array}
\end{eqnarray*}
and
\begin{itemize}

\item if $X_n\stackrel{P}{\longrightarrow} X$, then there exists a subsequence $\{X_{n_k}\}$ of $\{X_n\}$ such that $X_{n_k}\stackrel{a.s.}{\longrightarrow} X$ as $k\to\infty$;

\item if $X_n\stackrel{d}{\longrightarrow} C$, where $C$ is a constant, then $X_n\stackrel{P}{\longrightarrow} C$;

\item if $X_n\stackrel{d}{\longrightarrow} X$, then by  Skorokhod's theorem, there exist a sequence of random variables $\{Y,Y_n,n\geq 1\}$ such that for any $n\geq 1$, $X_n$ and $Y_n$ have the same distribution, $X$ and $Y$ have the same distribution, and $Y_n\stackrel{a.s.}{\longrightarrow} Y$.
\end{itemize}

Recently, motivated by the risk measures,
 superhedge pricing and modeling uncertain in finance, Peng \cite{Pe04}-\cite{Pe10}
 initiated the notion of independent and identically distributed (IID) random
 variables under sublinear expectations,  proved the weak law of large
 numbers and the central limit theorems, defined the $G$-expectations, $G$-Brownian motions and built It\^{o}'s type stochastic calculus.

In this note, we will survey the existing  convergence results for random variables under sublinear expectations, and prove some new results. Concretely, under the assumption that the sublinear expectation has the monotone continuity property, we will prove that $L^p$ convergence is stronger than convergence in capacity, convergence in capacity is stronger than convergence in distribution, and give some equivalent characterizations of convergence in distribution. In addition, a dominated convergence theorem under  sublinear expectations is given, which may have its own interest.

\section{Sublinear expectations}

In this section, we present some basic settings about sublinear expectations. Please refer to  Peng \cite{Pe04}-\cite{Pe10} for more details.

Let $\Omega$ be a given set and let $\bf{\mathcal {H}}$ be a linear space of real functions defined on $\Omega$ such that
 for any constant number $c, c\in \mathcal{H}$;  if $X\in\mathcal{H}$, then $|X|\in\mathcal{H}$; if $X_1,\ldots,X_n\in\mathcal{H}$, then for any $\varphi\in C_{l,Lip}(\mathbb{R}^n)$, $\varphi(X_1,\ldots,X_n)\in\mathcal{H}$, where $C_{l,Lip}(\mathbb{R}^n)$ denotes the linear space of functions $\varphi$ satisfying
$$|\varphi(x)-\varphi(y)|\le C(1+|x|^m+|y|^m)|x-y|, \forall x,y\in \mathbb{R}^n,$$
for some $C>0, m\in \mathbb{N}$ depending on $\varphi$.

\begin{defi}
A  sublinear expectation ${\hat{E}}$ on $\mathcal{H}$ is a functional $\hat{E}:\mathcal{H}\rightarrow\mathbb{R}$ satisfying the
following properties:\\
\hspace*{0.3cm} (a) Monotonicity: $\hat{E}[X]\ge\hat{E}[Y]$, if $X\ge Y$.\\
\hspace*{0.3cm}  (b) Constant preserving: $\hat{E}[c]=c,\forall c\in\mathbb{R}.$\\
\hspace*{0.3cm}  (c) Sub-additivity: $\hat{E}[X+Y]\le\hat{E}[X]+\hat{E}[Y].$\\
\hspace*{0.3cm} (d) Positive homogeneity: $\hat{E}[\lambda X]=\lambda\hat{E}[X]$, $\forall\lambda\ge 0.$\\
The triple $(\Omega,\mathcal{H},\hat{E})$ is called a
sublinear expectation space.
\end{defi}


For simplicity, we assume that $\Omega$ is a complete separable metric space, and use $\mathcal{B}(\Omega)$ to denote the Borel $\sigma$-algebra of $\Omega$. Further, we assume that there exists a family $\mathcal{P}$ of probability measures on $(\Omega,\mathcal{B}(\Omega))$ such that
$$
\hat{E}[X]=\sup_{P\in \mathcal{P}}E_P[X],\ \forall X\in \mathcal{H}.
$$
Suppose that for any $A\in \mathcal{B}(\Omega)$, $I_{A}\in\mathcal{H}$. Define
$$
V(A):={\hat{E}}[I_{A}]=\sup_{P\in\mathcal{P}}P(A), \forall A \in \mathcal{B}(\Omega).
$$
Obviously, $V(\emptyset)=0, V(\Omega)=1.$

\begin{thm} (\cite[Theorem 1]{DHP11} or \cite[Theorem VI.1.1]{Pe10}) The set function $V(\cdot)$ is a Choquet capacity, i.e.
\begin{itemize}
\item[1.] $0\leq V(A)\leq 1, \forall A\in \mathcal{B}(\Omega)$;

\item[2.] If $A\subset B$, then $V(A)\leq V(B)$;

\item[3.] If $(A_n)_{n=1}^{\infty}$ is a sequence in $\mathcal{B}(\Omega)$, then $V(\cup A_n)\leq \sum_n V(A_n)$;

\item[4.] If $(A_n)_{n=1}^{\infty}$ is an increasing sequence in $\mathcal{B}(\Omega): A_n\uparrow A=\cup A_n$, then $V(A)=\lim_{n\to\infty}V(A_n)$.
\end{itemize}
\end{thm}

\begin{defi} (\cite[Definition VI.1.3]{Pe10}) We use the standard capacity-related vocabulary: a set $A$ is {\bf polar} if $V(A)=0$ and a property holds ``{\bf quasi-surely}" (q.s.) if it holds outside a polar set.
\end{defi}

%
%

\begin{defi}(\cite[Definition 3.1]{CJP11})
For $p \in [1, \infty)$, the map
$$
{|| \cdot ||}_p : X \mapsto (\hat{E}[|X|^p])^{1/p}
$$
forms a seminorm on $\cal{H}$. Define the space ${\cal{L}}^{p}(\cal{F})$ as the completion under ${||\cdot||}_{p}$ of the set
$$
\{X\in \cal{H}: {||X||}_{p} < \infty\}
$$
and then $L^{p}(\cal{F})$ as the equivalence classes of ${\cal{L}}^{p}$ modulo equality in ${||\cdot||}_{p}$.
\end{defi}

\begin{defi}\label{defi2.4}(\cite[Definition 3.2]{CJP11})
Consider $K\subset L^{1}$. $K$ is said to be uniformly integrable
if $\hat{E}[I_{\{|X|\geq c\}}|X|]$ converges to 0 uniformly in $X\in K$ as
$c\rightarrow \infty$.
\end{defi}

\begin{defi}(\cite[Definition 3.3]{CJP11})
Let $L_{b}^{p}$ be the completion of the set of bounded functions $X \in \cal{H}$, under the norm ${|| \cdot ||}_p$ .
Note that $L_{b}^{p}\subset L^{p}$ .
\end{defi}

\begin{lem}(\cite[Proposition 15]{DHP11} or \cite[Lemma 3.4]{CJP11})\label{lem2.13}
For each $p\geq 1$,
$$
L_{b}^{p}=\{X\in L^{p}: \lim_{n\rightarrow \infty}\hat{E}[|X|^{p}I_{\{|X|>n\}}]=0\}.
$$
\end{lem}

\begin{thm}\label{thm2.9}
{\rm (Monotone Convergence Theorem)} (\cite[Theorem 2.2]{CJP11})
Let $\{X_n,n\geq 1\}$ be a sequence in $\mathcal{H}$ and lower bounded. If $X_n\uparrow X\in \mathcal{H}$. Then $\hat{E}[X_n]\uparrow \hat{E}[X]$.
\end{thm}

\section{Convergences under sublinear expectations}\setcounter{equation}{0}

Let $(\Omega,\mathcal{H},\hat{E})$ be a sublinear space introduced in Section 2. We further suppose that  for any $X\in \mathcal{H}$ and any bounded continuous function $f$, $f(X)\in \mathcal{H}$. In this section, we consider the convergences of random variables under sublinear expectations.
Let $\{X,X_n,n\geq 1\}$ be a sequence of random varialbes in $\mathcal{H}$. We have the following convergences:

\begin{itemize}

\item $\{X_n,n\geq 1\}$  is said to quasi-surely converge to $X$, if there exists a set $N\subset \Omega$ such that $\hat{E}[I_{N}]=0$ and  $\forall\omega\in \Omega\backslash N, \lim_{n\to\infty}X_n(\omega)=X(\omega)$, which is denoted by $X_n\stackrel{q.s.}{\longrightarrow} X$ or $X_n\to X\ q.s.$

\item $\{X_n,n\geq 1\}$  is said to converge to $X$ in capacity, if for any $\varepsilon>0$, $\lim_{n\to\infty}V(\{|X_n-X|\geq \varepsilon\})=0$,  which is denoted by $X_n\stackrel{V}{\longrightarrow} X$.

\item $\{X_n,n\geq 1\}$  is said to $L^p$ converge to $X$ $(p>0)$ if $\lim_{n\to\infty}\hat{E}[|X_n-X|^p]=0$, which is denoted by $X_n\stackrel{L^p}{\longrightarrow} X$.

\item $\{X_n,n\geq 1\}$  is said to converge to $X$ in distribution, if for any bounded continuous function $\varphi$, $\lim_{n\to\infty}\hat{E}[\varphi(X_n)]=\hat{E}[\varphi(X)]$,  which is denoted by $X_n\stackrel{d}{\longrightarrow} X$.

\item $\{X_n,n\geq 1\}$  is said to completely converge to $X$, if for any $\varepsilon>0$, $\sum_{n=1}^{\infty}V(\{|X_n-X|\geq \varepsilon\})<\infty$,  which is denoted by $X_n\stackrel{c.c.}{\longrightarrow} X$ .

\item $\{X_n,n\geq 1\}$  is said to s-$L^r$ converge to $X$ $(r>0)$ if $\sum_{n=1}^{\infty}\hat{E}[|X_n-X|^r]<\infty$, which is denoted by $X_n\stackrel{s\mbox{-}L^r}{\longrightarrow} X$ .
\end{itemize}

\subsection{Without the monotone continuity property}

(a) By  Markov's inequality (\cite[Lemma 10]{DHP11} or \cite[Lemma VI.1.13]{Pe10}), we get
$$
X_n\stackrel{s\mbox{-}L^r}{\longrightarrow} X \Rightarrow  X_n\stackrel{c.c.}{\longrightarrow} X.
$$

(b) By  the Borel-Cantell lemma (\cite[Lemma VI.1.5]{Pe10}), we get
$$
X_n\stackrel{c.c.}{\longrightarrow} X \Rightarrow  X_n\stackrel{q.s.}{\longrightarrow} X.
$$

(c) By  Markov's inequality (\cite[Lemma 10]{DHP11} or \cite[Lemma VI.1.13]{Pe10}), we get
$$
X_n\stackrel{L^p}{\longrightarrow} X \Rightarrow  X_n\stackrel{V}{\longrightarrow} X.
$$

(d) By \cite[Proposition VI.1.17]{Pe10}, we know that if  $X_n\stackrel{L^p}{\longrightarrow} X$, then there exists a subsequence $X_{n_k}$ such that $X_{n_k}\stackrel{q.s.}{\longrightarrow} X$.

In general, we do not have  ``$
X_n\stackrel{q.s.}{\longrightarrow} X \Rightarrow  X_n\stackrel{V}{\longrightarrow} X
$" and the dominated convergence theorem  as in the classic probability space.

By using the idea in the sufficiency proof of \cite[Theorem 3.2]{CJP11}, we  give the following  dominated convergence  theorem.

\begin{thm} \label{thm3.1}
Suppose that $\{X_n\}$ is a sequence in $L^1$, and $|X_n|\leq Y, \forall n\geq 1$ with $Y\in L_{b}^{1}$ . Further suppose that  $\{X_n\}$ converges to $X$ in capacity and $X\in L_b^1$. Then $\{X_n\}$ converges  to $X$ in $L^1$ norm.
\end{thm}
{\bf Proof.} The idea comes from the sufficiency proof of \cite[Theorem 3.2]{CJP11}. For the reader's convenience, we spell out the details.

By  Definition \ref{defi2.4} and Lemma \ref{lem2.13}, and the assumption that  $|X_n|\leq Y$ for all $n\geq 1$ and $Y\in
L_{b}^{1}$,   we get that $X_n\in L_{b}^{1}$ for all $n\geq 1$ and
$\{X_n\}$ is uniformly integrable.

For any $c>0$, we have
\begin{eqnarray}\label{thm3.1-a}
\hat{E}[|X_n-X|]&=&\hat{E}\left[\left|X_n(I_{\{|X_n|\leq c\}}+I_{\{|X_n|> c\}})-X(I_{\{|X|\leq c\}}+I_{\{|X|> c\}})\right|\right]\nonumber\\
&\leq& \hat{E}\left[\left|X_nI_{\{|X_n|\leq c\}}-XI_{\{|X|\leq c\}}\right|\right]+\hat{E}\left[|X_n|I_{\{|X_n|> c\}}\right]+
\hat{E}\left[|X|I_{\{|X|> c\}}\right]\nonumber\\
&\leq &\hat{E}\left[|X_n-X|I_{\{|X_n|\leq c,|X|\leq c\}}\right]+\hat{E}\left[|X_n|I_{\{|X_n|\leq c,|X|> c\}}\right]+\hat{E}\left[|X|I_{\{|X_n|> c,|X|\leq c\}}\right]\nonumber\\
&&+\hat{E}\left[|X_n|I_{\{|X_n|> c\}}\right]+\hat{E}\left[|X|I_{\{|X|> c\}}\right]\nonumber\\
&\leq&\hat{E}\left[|X_n-X|I_{\{|X_n|\leq c,|X|\leq c\}}\right]+\hat{E}\left[cI_{\{|X|> c\}}\right]+\hat{E}\left[cI_{\{|X_n|> c\}}\right]\nonumber\\
&&+\hat{E}\left[|X_n|I_{\{|X_n|> c\}}\right]+\hat{E}\left[|X|I_{\{|X|> c\}}\right]\nonumber\\
&\leq&\hat{E}\left[|X_n-X|I_{\{|X_n|\leq c,|X|\leq c\}}\right]+2\hat{E}\left[|X_n|I_{\{|X_n|> c\}}\right]+2\hat{E}\left[|X|I_{\{|X|> c\}}\right].
\end{eqnarray}
For any fixed $\epsilon>0$, as  $\{X_n\}$ is uniformly integrable and  $X\in L_b^1$, there exists $c>0$ such that
\begin{eqnarray}\label{thm3.1-b}
\sup_{n\geq 1}\hat{E}\left[|X_n|I_{\{|X_n|> c\}}\right]<\frac{\epsilon}{8},\ \ \hat{E}\left[|X|I_{\{|X|> c\}}\right]<\frac{\epsilon}{8}.
\end{eqnarray}
Since $\{X_n\}$ converges to $X$ in capacity, there exists $N$ such that for any $n\geq N$, $\hat{E}[I_{\{|X_n-X|>\epsilon/4\}}]<\epsilon/8c$. Obviously, $|X_n-X|I_{\{|X_n|\leq c,|X|\leq c\}}\leq 2c$. Hence,
for any $n\geq N$, we have
\begin{eqnarray}\label{thm3.1-c}
\hat{E}\left[|X_n-X|I_{\{|X_n|\leq c,|X|\leq c\}}\right]&\leq& \hat{E}\left[(|X_n-X|I_{\{|X_n|\leq c,|X|\leq c\}})\vee \frac{\epsilon}{4}\right]\nonumber\\
&\leq&\frac{\epsilon}{4}+2c\hat{E}\left[I_{\{|X_n-X|>\epsilon/4\}}\right]\leq \frac{\epsilon}{2}.
\end{eqnarray}
By (\ref{thm3.1-a})-(\ref{thm3.1-c}), we get that for any $n\geq N$, $\hat{E}[|X_n-X|]<\epsilon$. Hence $\{X_n\}$ converges  to $X$ in $L^1$ norm.\hfill\fbox

\begin{rem}\label{rem3.2}
If there exists a subsequence $\{X_{n_k}\}$ of $\{X_n\}$ such that $X_{n_k} \rightarrow X$ q.s. and $\{X_{n_k}\}$ is uniformly integrable, then $X\in L_b^1$. In fact,  by  Theorem \ref{thm2.9} and the monotonicity of $\hat{E}$, we get that for any $A\in \cal{F}$,
\begin{eqnarray}\label{rem3.2-a}
\hat{E}[I_A|X|]&=&\hat{E}[\lim_{k\to\infty}I_A|X_{n_k}|]=\hat{E}[\liminf_{k\to\infty}I_A|X_{n_k}|]\nonumber\\
&=&\hat{E}[\lim_{n\to\infty}\inf_{k\geq n}I_A|X_{n_k}|]\nonumber\\
&=&\lim_{n\to\infty}\hat{E}[\inf_{k\geq n}I_A|X_{n_k}|]\nonumber\\
&\leq&\lim_{n\to\infty}\inf_{k\geq n}\hat{E}[I_A|X_{n_k}|].
\end{eqnarray}
By  (\ref{rem3.2-a}) (setting $A=\Omega$), the fact that $\{X_{n_k}\}$ is uniformly integrable and  \cite[Theorem 3.1]{CJP11}, we get that
$\hat{E}[|X|]<\infty$. Then by the Markov's inequality, we get that
\begin{eqnarray}\label{rem3.2-b}
\lim_{m\to\infty}\hat{E}[I_{\{|X|>m\}}]\leq \lim_{m\to\infty}\frac{\hat{E}[|X|]}{m}=0.
\end{eqnarray}
By (\ref{rem3.2-a}), (\ref{rem3.2-b}), the fact that $\{X_{n_k}\}$ is uniformly integrable and  \cite[Theorem 3.1]{CJP11}, we get that
$$
\lim_{m\to\infty}\hat{E}[|X|I_{\{|X|>m\}}]=0,
$$
which together with Lemma \ref{lem2.13} implies that $X\in L_b^1$.
\end{rem}

\subsection{With the monotone continuity property}
Throughout the rest of the paper, we assume  that $\hat{E}$ has the monotone continuity property (\cite[Definition 2.2(vii)]{CJP11}), i.e. for any $X_n\downarrow 0$ on $\Omega$, we have $\hat{E}[X_n]\downarrow 0$. Then we have

\begin{rem}
$V(\cdot)$ is a continuous capacity, i.e. if $(A_n)_{n=1}^{\infty}$ is a decreasing sequence in $\mathcal{B}(\Omega): A_n\downarrow A=\cap A_n$, then $V(A)=\lim_{n\to\infty}V(A_n)$. In fact, we have $I_{A_n}-I_A=I_{A_n\backslash A}\downarrow 0$. Thus for any $n\geq 1$, we have
\begin{eqnarray*}
0\leq V(A_n)-V(A)=\hat{E}[I_{A_n}]-\hat{E}[I_{A}]\leq \hat{E}[I_{A_n}-I_{A}]\downarrow 0,
\end{eqnarray*}
which implies that $V(A)=\lim_{n\to\infty}V(A_n)$.
\end{rem}

By \cite[Lemma 3.7]{CJP11}, we know that $
X_n\stackrel{q.s.}{\longrightarrow} X \Rightarrow  X_n\stackrel{V}{\longrightarrow} X,
$ and if $X_n\stackrel{V}{\longrightarrow} X$, then there exists a subsequence $X_{n_k}$ such that $X_{n_k}\stackrel{q.s.}{\longrightarrow} X$. Thus, by Remark \ref{rem3.2} we can rewrite Theorem \ref{thm3.1} in this case as follows.

\begin{thm}\label{dcthm} (\cite[Theorem 3.11]{Chen12}) \label{thm2.16}
Suppose that $\{X_n\}$ is a sequence in $L^1$, and $|X_n|\leq Y, \forall n\geq 1$ with $Y\in L_{b}^{1}$ . If $X_n \rightarrow X$ q.s.  or $\{X_n\}$ converges to $X$ in capacity, then $\{X_n\}$ converges  to $X$ in $L^1$ norm.
\end{thm}

\subsubsection{Convergence in capacity}
In this subsection, we give some discussions of convergence in capacity.

\begin{pro}\label{lem3.1}
$X_n\stackrel{V}{\longrightarrow} X$ if and only  if for any subsequence $\{X_{n'}\}$ of $\{X_{n}\}$, there exists a subsequence $\{X_{n'_k}\}$ of $\{X_{n'}\}$ such that $X_{n'_k}\stackrel{q.s.}{\longrightarrow} X$.
\end{pro}
{\bf Proof.} (Necessity) Suppose that $\{X_{n'}\}$ is a subsequence of $\{X_n\}$. Then $X_n'\stackrel{V}{\longrightarrow} X$, and thus by \cite[Lemma 3.7]{CJP11}, there exists a subsequence $\{X_{n'_k}\}$ of $\{X_{n'}\}$ such that $X_{n'_k}\stackrel{q.s.}{\longrightarrow} X$.

(Sufficiency) We use proof by contradiction. We assume that $\{X_n\}$ does not converge to $X$ in capacity. Then by the definition, there exist $\varepsilon>0$ and $\delta>0$ such that
$$
\limsup_{n\to\infty}V(\{|X_n-X|\geq \varepsilon\})>\delta.
$$
It follows that there exists a subsequence $\{X_{n_k}\}$ such that for any $k\geq 1$, $V(\{|X_{n_k}-X|\geq \varepsilon\})>\delta$. By the assumption, the exists a subsequence $\{X_{n_{k_l}}\}$ of $\{X_{n_k}\}$ such that $X_{n_{k_l}}\stackrel{q.s.}{\longrightarrow} X$, which implies that
$$
V(\cap_{n=1}^{\infty}\cup_{l=n}^{\infty}\{|X_{n_{k_l}}-X|\geq \varepsilon\})=0,
$$
i.e.
\begin{eqnarray}\label{lem3.1-a}
\lim_{n\to\infty}V(\cup_{l=n}^{\infty}\{|X_{n_{k_l}}-X|\geq \varepsilon\})=0.
\end{eqnarray}
But for any $m\geq n$, we have
$$
V(\cup_{l=n}^{\infty}\{|X_{n_{k_l}}-X|\geq \varepsilon\})\geq V(\{|X_{n_{k_m}}-X|\geq \varepsilon\})>\delta>0.
$$
It contradicts with (\ref{lem3.1-a}). Hence we obtain that $\{X_n\}$  converges to $X$ in capacity.\hfill\fbox

\begin{cor}\label{cor3.3}
  Suppose that $X_n\stackrel{V}{\longrightarrow} X$, and $f$ is a continuous function such that $\{f(X),f(X_n),n\geq 1\}\subset \mathcal{H}$. Then  $f(X_n)\stackrel{V}{\longrightarrow} f(X)$.
\end{cor}
{\bf Proof.} Suppose that $\{f(X_{n_k})\}$ is a subsequence of $\{f(X_n)\}$. Then by the assumption and Proposition \ref{lem3.1}, there exists a subsequence $\{X_{n_{k_l}}\}$ of $\{X_{n_k}\}$ such that $X_{n_{k_l}}\stackrel{q.s.}{\longrightarrow} X$. Since $f$ is a continuous function, we get that $f(X_{n_{k_l}})\stackrel{q.s.}{\longrightarrow} f(X)$. By Proposition \ref{lem3.1} again, we obtain that $f(X_n)\stackrel{V}{\longrightarrow} f(X)$.\hfill\fbox

\subsubsection{Convergence in distribution}

In this subsection, we discuss  convergence in distribution. At first, we prove that convergence in capacity is stronger than convergence in distribution.

\begin{pro}\label{pro3.7}
 $X_n\stackrel{V}{\longrightarrow} X \Rightarrow  X_n\stackrel{d}{\longrightarrow} X.$
\end{pro}
{\bf Proof.} Suppose that $X_n\stackrel{V}{\longrightarrow} X$ and $f$ is a bounded continuous function. By Corollary \ref{cor3.3}, we get that $f(X_n)\stackrel{V}{\longrightarrow} f(X)$. Then by the dominated convergence theorem (see Theorem \ref{dcthm} above), we obtain that
$$
\lim_{n\to\infty}\hat{E}[|f(X_n)-f(X)|]=0,
$$
which implies
$$
\lim_{n\to\infty}\hat{E}[f(X_n)]=\hat{E}[f(X)].
$$
Hence $X_n\stackrel{d}{\longrightarrow} X$. \hfill\fbox

\bigskip

Now we have the following relations under sublinear expectations as in the classic setting:
\begin{eqnarray*}
\begin{array}{ccccccccc}
X_n\stackrel{s\mbox{-}L^r}{\longrightarrow} X &\Rightarrow & X_n\stackrel{c.c.}{\longrightarrow} X& \Rightarrow & X_n\stackrel{q.s.}{\longrightarrow} X &\Rightarrow & X_n\stackrel{V}{\longrightarrow} X &\Rightarrow & X_n\stackrel{d}{\longrightarrow} X.\\
&&&&&&\Uparrow&&\\
&&&&&&X_n\stackrel{L^p}{\longrightarrow} X&&
\end{array}
\end{eqnarray*}

 In order to state next result, we introduce one notion and one assumption.

\begin{defi}\label{defi3.5}
Let $X$ be a random variable in $\mathcal{H}$. Define a set function $C_X$ on $(\mathbf{R},\mathcal{B}(\mathbf{R}))$ as follows:
$$
C_X(A)=\hat{E}[I_{X^{-1}(A)}],\ \forall A\in \mathcal{B}(\mathbf{R}).
$$
Then $C_X$ is a continuous capacity, and we call it the distribution capacity of $X$.
\end{defi}

 For simplicity, we define

{\it Assumption (A)}: The set $\{x\in \mathbf{R}: \hat{E}[I_{\{X=x\}}]>0\}$ is at most countable.

\begin{rem}\label{rem3.6}
Let $X$ be a random variable in $\mathcal{H}$ and $C_X$ be its distribution capacity. Suppose that $C_X$  is 2-monotone, i.e.
$$
C_X(A\cup B)+C_X(A\cap B)\geq C_X(A)+C_X(B).
$$
It follows  that if $\{A_k,k=1,2,\ldots,n\}$ satisfies that $A_i \cap A_j=\emptyset,\forall i\neq j$, then
$$
C_X(\cup_{k=1}^nA_k)\geq \sum_{k=1}^nC_X(A_k).
$$
It follows that $X$ satisfies Assumption (A)  in this case.
\end{rem}

\begin{thm}\label{thm3.7}
Suppose that $\{X,X_n,n\geq 1\}$ is a sequence of random variables in the sublinear expectation space $(\Omega,\mathcal{H},\hat{E})$. Define the following six claims:
\begin{itemize}
\item[(1)] $X_n\stackrel{d}{\longrightarrow} X$;

\item[(2)] For any bounded lower semi-continuous function $f:\mathbf{R}\to \mathbf{R}, \liminf_{n\to\infty}\hat{E}[f(X_n)]\geq \hat{E}[f(X)]$, and for any bounded upper semi-continuous function $g:\mathbf{R}\to \mathbf{R}, \limsup_{n\to\infty}\hat{E}[g(X_n)]\leq \hat{E}[g(X)]$;

    \item[(3)] For any bounded $X$-q.s. continuous function $f: \mathbf{R}\to \mathbf{R}$ (i.e. if  $D_f:=\{x\in \mathbf{R}|f\ \mbox{is not continuous}\linebreak \mbox{at}\ x\}$, then $\{\omega\in \Omega|X(\omega)\in D_f\}$ is polar), $\lim_{n\to\infty}\hat{E}[f(X_n)]=\hat{E}[f(X)]$;

\item[(4)] For any open subset $A$ of $\Omega$,
$$
\liminf_{n\to\infty}\hat{E}[I_A(X_n)]\geq \hat{E}[I_A(X)],\  \limsup_{n\to\infty}\hat{E}[-I_A(X_n)]\leq \hat{E}[-I_A(X)];
 $$

\item[(5)] For any  closed subset $F$ of $\Omega$,
$$
\limsup_{n\to\infty}\hat{E}[I_F(X_n)]\leq \hat{E}[I_F(X)],\ \liminf_{n\to\infty}\hat{E}[-I_F(X_n)]\geq \hat{E}[-I_F(X)];
 $$

\item[(6)] For any $X$-q.s. continuous set $A$ (i.e. $\{\omega\in \Omega|X(\omega)\in \partial A\}$ is polar), $\hat{E}[I_A(X_n)]\to \hat{E}[I_A(X)]$ and $\hat{E}[-I_A(X_n)]\to \hat{E}[-I_A(X)]$.
\end{itemize}

Then we have
$$
(1)\Leftrightarrow (2) \Leftrightarrow (3) \Rightarrow (4)\Leftrightarrow  (5) \Rightarrow (6).
$$
Further, if Assumption (A) holds,  then (6) $\Rightarrow$ (1).  Hence, in this case, (1)-(6) are equivalent.
\end{thm}

\noindent{\bf Proof.} (1) $\Rightarrow$ (2): If $f'\leq f$ and $f'$ is a bounded continuous function, then by (1), we have
\begin{eqnarray}\label{thm3.7-a}
\liminf_{n\to\infty}\hat{E}[f(X_n)]\geq \liminf_{n\to\infty}\hat{E}[f'(X_n)]=\hat{E}[f'(X)].
\end{eqnarray}
Suppose that $f$ is a bounded lower semi-continuous function, then there exists a sequence $\{f_n\}$   of increasing bounded continuous functions such that $f_n\uparrow f$. By (\ref{thm3.7-a}) and the monotone convergence theorem (see Theorem \ref{thm2.9} above), we get that
$$
\liminf_{n\to\infty}\hat{E}[f(X_n)]\geq \lim_{m\to\infty}\hat{E}[f_m(X)]=\hat{E}[f(X)].
$$

If $g'\geq g$ and $g'$ is a bounded continuous function, then by (1), we have
\begin{eqnarray}\label{thm3.7-b}
\limsup_{n\to\infty}\hat{E}[g(X_n)]\leq \limsup_{n\to\infty}\hat{E}[g'(X_n)]=\hat{E}[g'(X)].
\end{eqnarray}
Suppose that $g$ is a bounded upper semi-continuous function, then there exists a sequence $\{g_n\}$   of decreasing bounded continuous functions such that $g_n\downarrow g$. By (\ref{thm3.7-b}), we get that
\begin{eqnarray}\label{thm3.7-c}
\limsup_{n\to\infty}\hat{E}[g(X_n)]\leq \lim_{m\to\infty}\hat{E}[g_m(X)].
\end{eqnarray}
Since $g_m\downarrow g$, we have $g_m(X)-g(X)\downarrow 0$, and thus
$$
0\leq \hat{E}[g_m(X)]-\hat{E}[g(X)]\leq \hat{E}[g_m(X)-g(X)]\downarrow 0.
$$
It follows that
$$
\lim_{m\to\infty}\hat{E}[g_m(X)]=\hat{E}[g(X)],
$$
which together with (\ref{thm3.7-c}) implies that
\begin{eqnarray*}
\limsup_{n\to\infty}\hat{E}[g(X_n)]\leq \hat{E}[g(X)].
\end{eqnarray*}
Hence (2) holds.

(2) $\Rightarrow$ (3): Suppose that $f$ is a function in (3). Define
\begin{eqnarray*}
&&\underline{f}(x)=\liminf_{y\to x}f(y)=\sup_{n\geq 1}\inf_{y\in B(x,1/n)}f(y),\\
&&\overline{f}(x)=\limsup_{y\to x}f(y)=\inf_{n\geq 1}\sup_{y\in B(x,1/n)}f(y),
\end{eqnarray*}
where $B(x,1/n)$ stands for the open ball $\{y\in \mathbf{R}: |y-x|<1/n\}$. Then $\underline{f}$ is a bounded lower semi-continuous function,  $\overline{f}$ is a bounded upper semi-continuous function, $\underline{f}(x)\leq f(x)\leq \overline{f}(x),\forall x$, and if $x$ is a continuous point of $f$, $\underline{f}(x)=f(x)= \overline{f}(x)$. By the assumption on $f$, we have
$$
\hat{E}[\underline{f}(X)]=\hat{E}[f(X)]=\hat{E}[\overline{f}(X)],
$$
which together with (2) implies that
\begin{eqnarray*}
\hat{E}[f(X)]&=&\hat{E}[\underline{f}(X)]\leq \liminf_{n\to\infty}\hat{E}[\underline{f}(X_n)]\\
&\leq& \liminf_{n\to\infty}\hat{E}[f(X_n)]\leq \limsup_{n\to\infty}\hat{E}[f(X_n)]\\
&\leq&\limsup_{n\to\infty}\hat{E}[\overline{f}(X_n)]\\
&\leq& \hat{E}[\overline{f}(X)]=\hat{E}[f(X)].
\end{eqnarray*}
Hence (3) holds.

(3) $\Rightarrow$ (1): It is obvious.

(3) $\Rightarrow$ (4):  Suppose that (3) holds. Then (2) holds.  If $A$ is an open set, then $I_A$ is a bounded lower semi-continuous function, and  $-I_A$ is a bounded upper semi-continuous function. Then by (2), we obtain (4).

(4) $\Leftrightarrow$ (5): Suppose that (4) holds and $F$ is a closed set. Then $I_F=1-I_{F^c}$, where $F^c=\Omega-F$ is an open set. Further, we have
\begin{eqnarray*}
&&\hat{E}[I_F(X_n)]=1+\hat{E}[-I_{F^c}(X_n)],\ \hat{E}[I_F(X)]=1+\hat{E}[-I_{F^c}(X)];\\
&&\hat{E}[-I_F(X_n)]=\hat{E}[I_{F^c}(X_n)]-1,\ \hat{E}[-I_F(X)]=\hat{E}[I_{F^c}(X)]-1.
\end{eqnarray*}
It follows that (5) holds. Similarly, we can prove (5) $\Rightarrow$ (4).

 (5) $\Rightarrow$ (6): Suppose that (5) holds. Then (4) holds also.   Denote by $A^o$ the interior of $A$, by $\bar{A}$ the closure of $A$. Then $A^o\subset A\subset \bar{A}$, $A^o$ is an open set and $\bar{A}$ is a closed set. By the assumption on $A$, we have
$$
\hat{E}[I_{A^o}(X)]=\hat{E}[I_A(X)]=\hat{E}[I_{\bar{A}}(X)],
$$
which together with (4) and (5) implies that
$$
\limsup_{n\to\infty}\hat{E}[I_{A}(X_n)]\leq \limsup_{n\to\infty}\hat{E}[I_{\bar{A}}(X_n)]\leq \hat{E}[I_{\bar{A}}(X)]=\hat{E}[I_A(X)],
$$
and
$$
\liminf_{n\to\infty}\hat{E}[I_{A}(X_n)]\geq \limsup_{n\to\infty}\hat{E}[I_{A^o}(X_n)]\geq \hat{E}[I_{A^o}(X)]=\hat{E}[I_A(X)].
$$
Hence $\hat{E}[I_A(X_n)]\to \hat{E}[I_A(X)]$. Similarly, we can get $\hat{E}[-I_A(X_n)]\to \hat{E}[-I_A(X)]$. Hence (6) holds.

(6) $\Rightarrow$ (1)  under Assumption (A):  Let $f\in C_{b}(\mathbf{R})$.  Then there exists $M>0$ such that $|f|<M$. Define $B=\{c\in \mathbf{R}|C_X(f^{-1}(c))\neq 0\}$. By  Assumption (A), we get that $B$ is at most countable.  For any $\varepsilon>0$, we construct a partition $\{t_i\}_{i=1}^{j}$ of $[-M,M]$  satisfying that $
-M=t_{0}<t_{1}<\cdot\cdot\cdot<t_{j}=M,t_i\notin B,$ and $\sup_{i}(t_{i}-t_{i-1})<\varepsilon,i=1,2,\ldots,j$. Define $B_i=\{x:t_i\leq f(x)<t_{i+1}\},i=0,1,\ldots,j-1,$ then $B_{i}$ is an $X$-q.s. continuous set, and for any $X\in \mathcal{H}$,$|f(X)-\sum_{i=0}^{j-1}t_iI_{B_i}(X)|<\varepsilon.$

By (6), we know that
$$
\hat{E}[I_{B_i}(X_n)]\to \hat{E}[I_{B_i}(X)],\ \hat{E}[-I_{B_i}(X_n)]\to \hat{E}[-I_{B_i}(X)].
$$
It follows that for any constant $\alpha$,   $
\hat{E}[\alpha I_{B_i}(X_n)]\to \hat{E}[\alpha I_{B_i}(X)],$
and thus
\begin{eqnarray}\label{thm3.7-d}
\sum_{i=0}^{j-1}\hat{E}[t_iI_{B_i}(X_n)]\to \sum_{i=0}^{j-1}\hat{E}[t_iI_{B_i}(X)].
\end{eqnarray}

By the triangle inequality  and  the sublinear property, we have
\begin{eqnarray*}
|\hat{E}[f(X_n)]-\hat{E}[f(X)]|&\leq& \left|\hat{E}[f(X_n)]-\sum_{i=0}^{j-1}\hat{E}[t_iI_{B_i}(X_n)]\right|+\left|\sum_{i=0}^{j-1}\hat{E}[t_iI_{B_i}(X_n)]- \sum_{i=0}^{j-1}\hat{E}[t_iI_{B_i}(X)]\right|\\
&&+\left|\hat{E}[f(X)]-\sum_{i=0}^{j-1}\hat{E}[t_iI_{B_i}(X)]\right|\\
&\leq&\hat{E}\left[\left|\left(f-\sum_{i=0}^{j-1}t_iI_{B_i}\right)(X_n)\right|\right]+\left|\sum_{i=0}^{j-1}\hat{E}[t_iI_{B_i}(X_n)]- \sum_{i=0}^{j-1}\hat{E}[t_iI_{B_i}(X)]\right|\\
&&+\hat{E}\left[\left|\left(f-\sum_{i=0}^{j-1}t_iI_{B_i}\right)(X)\right|\right]\\
&\leq& 2\varepsilon+\left|\sum_{i=0}^{j-1}\hat{E}[t_iI_{B_i}(X_n)]- \sum_{i=0}^{j-1}\hat{E}[t_iI_{B_i}(X)]\right|,
\end{eqnarray*}
which together with (\ref{thm3.7-d}) implies that
$$
\limsup_{n\to\infty}|\hat{E}[f(X_n)]-\hat{E}[f(X)]|\leq 2\varepsilon.
$$
By the arbitrariness of $\varepsilon$, we get
$$
\lim_{n\to\infty}|\hat{E}[f(X_n)]-\hat{E}[f(X)]|=0.
$$
 Hence (1) holds. \hfill\fbox

\begin{defi}\label{defi3.8}
Let $X$ be a random variable  in the sublinear expectation space $(\Omega,\mathcal{H},\hat{E})$. For any $x\in \mathbf{R}$, define
$$
\overline{F}(x)=\hat{E}[I_{\{X\leq x\}}],\ \underline{F}(x)=-\hat{E}[-I_{\{X\leq x\}}],
$$
where $\{X\leq x\}:=\{\omega\in\Omega|X(\omega)\leq x\}$. We call $(\overline{F},\underline{F})$ the distribution function pair of $X$.
\end{defi}

\begin{pro}\label{pro3.9}
Let $(\overline{F},\underline{F})$ be the distribution function pair of $X$. Then\\
(i) $0\leq \underline{F}(x)\leq \overline{F}(x)\leq 1,\ \forall x\in \mathbf{R}$;\\
(ii) both $\overline{F}$ and $\underline{F}$ are increasing functions and right continuous;\\
(iii) $\lim_{x\to-\infty}\overline{F}(x)=0,\ \lim_{x\to+\infty}\overline{F}(x)=1$;\ $\lim_{x\to-\infty}\underline{F}(x)=0,\ \lim_{x\to+\infty}\underline{F}(x)=1$;\\
(iv) if $\hat{E}[I_{\{X= x\}}]=0$, then $x$ is a continuous point of $\overline{F}$;\ if $\hat{E}[-I_{\{X= x\}}]=0$, then $x$ is a continuous point of $\underline{F}$.
\end{pro}
{\bf Proof.} (i) For any $x$, we have $0\leq I_{\{X\leq x\}}\leq 1$. Then  we have
$$
0\leq -\hat{E}[-I_{\{X\leq x\}}]\leq \hat{E}[I_{\{X\leq x\}}]\leq 1,
$$
i.e. (i) holds.

In the following, we only prove the results for $\overline{F}$ and the proofs for $\underline{F}$ are similar.

(ii) Obviously, $\overline{F}$ is an increasing function. When $y\downarrow x$, we have
$$
0\leq I_{\{X\leq y\}}-I_{\{X\leq x\}}\downarrow 0,
$$
which together with the monotone continuity property of $\hat{E}$ implies that
$$
0\leq \hat{E}[I_{\{X\leq y\}}]-\hat{E}[I_{\{X\leq x\}}]\leq \hat{E}[I_{\{X\leq y\}}-I_{\{X\leq x\}}]\downarrow 0,
$$
i.e.  $\lim_{y\downarrow x}\overline{F}(y)=\overline{F}(x)$. Hence $\overline{F}$ is right continuous.

(iii) When $x\to -\infty$, we have $I_{\{X\leq x\}}\downarrow 0$. Then by the monotone continuity property of $\hat{E}$, we get
$$
\lim_{x\to-\infty}\hat{E}[I_{\{X\leq x\}}]=0,
$$
i.e. $\lim_{x\to-\infty}\overline{F}(x)=0$. When $x\to +\infty$, we have $I_{\{X\leq x\}}\uparrow 1$. Then by the monotone convergence theorem (Theorem \ref{thm2.9}), we get
$$
\lim_{x\to+\infty}\hat{E}[I_{\{X\leq x\}}]=1,
$$
i.e. $\lim_{x\to+\infty}\overline{F}(x)=1$.

(iv) If $\hat{E}[I_{\{X= x\}}]=0$, then
$$
0\leq \hat{E}[I_{\{X\leq x\}}]-\hat{E}[I_{\{X<x\}}]\leq \hat{E}[I_{\{X\leq x\}}-I_{\{X<x\}}]=\hat{E}[I_{\{X= x\}}]=0.
$$
It follows that $\hat{E}[I_{\{X\leq x\}}]=\hat{E}[I_{\{X<x\}}]$, i.e. $\overline{F}(x)=\overline{F}(x-)$. Hence $x$ is a continuous point of $\overline{F}$. \hfill\fbox

%
%
%

\begin{thm}\label{thm3.10}
Suppose that $\{X,X_n,n\geq 1\}$ is a sequence of random variables in the sublinear expectation space $(\Omega,\mathcal{H},\hat{E})$ and  the distribution function pairs are $(\overline{F},\underline{F}),(\overline{F}_n,\underline{F}_n),n\geq 1$, respectively. Define the following three claims:
\begin{itemize}
\item[(1)] $X_n\stackrel{d}{\longrightarrow} X$;

\item[(2)] for a sequence $\{a_k,b_k,k=1,2,\ldots,n\}$ satisfying $a_1<b_1<a_2<b_2<\cdots<a_n<b_n$ and  $\hat{E}[I_{\{X= a_k\}}]=\hat{E}[I_{\{X= b_k\}}]=0,\forall k=1,2,\ldots,n$, it holds that $\lim_{n\to\infty}\hat{E}[I_{\cup_{k=1}^n (a_k,b_k]}(X_n)]=\hat{E}[I_{\cup_{k=1}^n (a_k,b_k]}(X)], \lim_{n\to\infty}\hat{E}[-I_{\cup_{k=1}^n (a_k,b_k]}(X_n)]=\hat{E}[-I_{\cup_{k=1}^n (a_k,b_k]}(X)]$;

\item[(3)] for any point $x$ with $\hat{E}[I_{\{X= x\}}]=0$, it holds that $\lim_{n\to\infty}\overline{F}_n(x)=\overline{F}(x)$ and $\lim_{n\to\infty}\underline{F}_n(x)=\underline{F}(x)$.
\end{itemize}
Then we have
$$
(1)\Rightarrow (2)\ \ \mbox{and}\ \ (1) \Rightarrow (3).
$$
 If $X$ satisfies Assumption (A), then (2) $\Rightarrow$ (1).
\end{thm}
{\bf Proof.} (1) $\Rightarrow$ (2): For  a sequence $\{a_k,b_k,k=1,2,\ldots,n\}$ satisfying the conditions in (2), the set $\cup_{k=1}^n (a_k,b_k]$ is an $X$-q.s. continuous set. Then by Theorem \ref{thm3.7}, we get  (2).

(1) $\Rightarrow$ (3): For any point $x$ with $\hat{E}[I_{\{X= x\}}]=0$, the set $(-\infty,x]$ is an $X$-q.s. continuous set. Then by Theorem \ref{thm3.7}, we get  (3).

(2) $\Rightarrow$ (1) under Assumption (A): Suppose that (2) and Assumption (A) hold. Let $A$ be an open subset of $\mathbf{R}$. Then $A$ can be expressed to be the disjoint union of open intervals $\{I_k,k\in I\}$. Without loss of generality, we assume that $I=\mathbb{N}$.

(a) For any $m\in \mathbb{N}$, we have
\begin{eqnarray}\label{thm3.10-a}
\liminf_{n\longrightarrow\infty}\hat{E}[I_{A}(X_{n})]&=&\liminf_{n\longrightarrow\infty}\hat{E}[I_{\cup_{k=1}^{\infty}I_k}(X_{n})]\nonumber\\
&\geq& \liminf_{n\longrightarrow\infty}\hat{E}[I_{\cup_{k=1}^mI_k}(X_{n})]\nonumber\\
&=&\liminf_{n\longrightarrow\infty}\hat{E}[\sum_{k=1}^mI_{I_k}(X_{n})].
\end{eqnarray}
Given $\varepsilon>0$. For any $k$, by Assumption (A) and the monotone continuity property of the sublinear expectation, there exists a  sub-interval $(a_k,b_k]$ of $I_k$ such that $\hat{E}[I_{\{X= a_k\}}]=\hat{E}[I_{\{X= b_k\}}]=0$, and
\begin{eqnarray}\label{thm3.10-a-1}
 \hat{E}[I_{I_k}(X)-I_{(a_k,b_k]}(X)]<\varepsilon/2^k.
\end{eqnarray}
It follows that for any $m\geq 1$, we have
\begin{eqnarray}\label{thm3.10-b}
0&\leq& \hat{E}[I_{\cup_{k=1}^mI_k}(X)]-\hat{E}[I_{\cup_{k=1}^m(a_k,b_k]}(X)]\nonumber\\
&\leq&\hat{E}[[I_{\cup_{k=1}^mI_k}(X)-I_{\cup_{k=1}^m(a_k,b_k]}(X)]\nonumber\\
&\leq &\sum_{k=1}^m\hat{E}[I_{I_k}(X)-I_{(a_k,b_k]}(X)]\nonumber\\
&\leq&\sum_{k=1}^m\varepsilon/2^k.
\end{eqnarray}
By (\ref{thm3.10-a}), the monotone property of $\hat{E}$,  the condition (2), (\ref{thm3.10-b}) and the monotone convergence theorem (Theorem \ref{thm2.9}), we obtain that
\begin{eqnarray*}
\liminf_{n\longrightarrow\infty}\hat{E}[I_{A}(X_{n})]
&\geq&\lim_{m\to\infty}\liminf_{n\longrightarrow\infty}\hat{E}[\sum_{k=1}^mI_{I_k}(X_{n})]\nonumber\\
&\geq&\lim_{m\to\infty}\liminf_{n\longrightarrow\infty}\hat{E}[\sum_{k=1}^mI_{(a_k,b_k]}(X_{n})]\nonumber\\
&=&\lim_{m\to\infty}\hat{E}[\sum_{k=1}^mI_{(a_k,b_k]}(X)]\nonumber\\
&\geq&\lim_{m\to\infty}\hat{E}[\sum_{k=1}^mI_{I_k}(X)]-\varepsilon\nonumber\\
&=&\hat{E}[I_A(X)]-\varepsilon.
\end{eqnarray*}
By the arbitrariness of $\varepsilon$, we get that
\begin{eqnarray}\label{thm3.10-c}
\liminf_{n\longrightarrow\infty}\hat{E}[I_{A}(X_{n})]\geq\hat{E}[I_A(X)].
\end{eqnarray}

(b) Given $\varepsilon>0$. For any $k$, take the  interval $(a_k,b_k]$ as in (a).  For any $m$, we have
\begin{eqnarray}\label{thm3.10-d}
\limsup_{n\to\infty}\hat{E}[-I_A(X_n)]
&\leq &\limsup_{n\to\infty}\hat{E}[-I_{\cup_{k=1}^mI_k}(X_n)]
\leq\limsup_{n\to\infty}\hat{E}[-I_{\cup_{k=1}^m(a_k,b_k]}(X_n)],
\end{eqnarray}
and by (\ref{thm3.10-a-1}), we have
\begin{eqnarray}\label{thm3.10-e}
0&\leq& \hat{E}[-I_{\cup_{k=1}^m(a_k,b_k]}(X)]-\hat{E}[-I_{\cup_{k=1}^mI_k}(X)]\nonumber\\
&\leq&\hat{E}[-I_{\cup_{k=1}^m(a_k,b_k]}(X)-(-I_{\cup_{k=1}^mI_k}(X))]\nonumber\\
&\leq &\sum_{k=1}^m\hat{E}[I_{I_k}(X)-I_{(a_k,b_k]}(X)]\nonumber\\
&\leq&\sum_{k=1}^m\varepsilon/2^k.
\end{eqnarray}
By the monotone continuity property of the sublinear expectation, we have
\begin{eqnarray}\label{thm3.10-f}
0&\leq& \hat{E}[-I_{\cup_{k=1}^mI_k}(X)]-\hat{E}[-I_A(X)]\nonumber\\
&\leq& \hat{E}[-I_{\cup_{k=1}^mI_k}(X)-(-I_A(X))]\nonumber\\
&=&\hat{E}[I_{\cup_{k=m+1}^{\infty}I_k}(X)\downarrow 0\ \mbox{as}\ m\to\infty.
\end{eqnarray}
By (\ref{thm3.10-d}), the condition (2), (\ref{thm3.10-e})  and (\ref{thm3.10-f}), we get
\begin{eqnarray*}
\limsup_{n\to\infty}\hat{E}[-I_A(X_n)]
&\leq &\lim_{m\to\infty}\limsup_{n\to\infty}\hat{E}[-I_{\cup_{k=1}^m(a_k,b_k]}(X_n)]\nonumber\\
&\leq&\lim_{m\to\infty}\hat{E}[-I_{\cup_{k=1}^m(a_k,b_k]}(X)]\nonumber\\
&\leq&\lim_{m\to\infty}\hat{E}[-I_{\cup_{k=1}^mI_k}(X)]+\varepsilon\nonumber\\
&=&\hat{E}[-I_A(X)]+\varepsilon.
\end{eqnarray*}
By the arbitrariness of $\varepsilon$, we get that
\begin{eqnarray}\label{thm3.10-g}
\limsup_{n\longrightarrow\infty}\hat{E}[-I_{A}(X_{n})]\leq\hat{E}[-I_A(X)].
\end{eqnarray}

By (\ref{thm3.10-c}), (\ref{thm3.10-g}) and Theorem \ref{thm3.7}, we obtain that (1) holds.
%
%
\hfill\fbox

\begin{pro}\label{pro3.12}
Suppose that $\{X,X_n,n\geq 1\}$ is a sequence of random variables in the sublinear expectation space $(\Omega,\mathcal{H},\hat{E})$ and $X_n\stackrel{d}{\longrightarrow} C$, where $C$ is a constant. Then $X_n\stackrel{V}{\longrightarrow} C$.
\end{pro}
{\bf Proof.} Denote by  $(\overline{F},\underline{F})$ the distribution function pairs  of $X$, and by $(\overline{F}_n,\underline{F}_n)$ the distribution function pairs of $X_n$ for any $n\geq 1$. We have
\begin{eqnarray*}
\overline{F}(x)=\underline{F}(x)=\left\{
\begin{array}{ll}
1, & \mbox{when}\ x\geq c,\\
0, &\mbox{otherwise},
\end{array}
\right.
\end{eqnarray*}
and for any $x\neq c$, we have $\hat{E}[I_{\{X=x\}}]=0.$ Then by Theorem \ref{thm3.10}, we get that for any $x\neq c$,
\begin{eqnarray}\label{pro3.12-a}
\lim_{n\to\infty}\overline{F}_n(x)=\lim_{n\to\infty}\underline{F}_n(x)=\left\{
\begin{array}{ll}
1, & \mbox{when}\ x> c,\\
0, &\mbox{when}\ x< c.
\end{array}
\right.
\end{eqnarray}

For any $\varepsilon>0$, we have
\begin{eqnarray*}
 \hat{E}[I_{\{|X_{n}-C|\geq \varepsilon\}}]&=&\hat{E}[I_{\{X_{n}\geq C+\varepsilon\}\bigcup\{X_{n}\leq C-\varepsilon\}}]\\
 & \leq& \hat{E}[I_{\{X_{n}\geq C+\varepsilon\}}]+\hat{E}[I_{\{X_{n}\leq C-\varepsilon\}}]\\
 &=&\hat{E}[1-I_{\{X_{n}<C+\varepsilon\}}]+\hat{E}[I_{\{X_{n}\leq C-\varepsilon\}}]\\
 &\leq& 1+\hat{E}[-I_{\{X_{n}< C+\varepsilon\}}]+\hat{E}[I_{\{X_{n}\leq C-\varepsilon\}}]\\
 &\leq& 1+\hat{E}[-I_{\{X_{n}\leq  C+\frac{\varepsilon}{2}\}}]+\hat{E}[I_{\{X_{n}\leq C-\varepsilon\}}]\\
 &=&1-\underline{F}_n(C+\frac{\varepsilon}{2})+\overline{F}_n(C-\varepsilon),
 \end{eqnarray*}
 which together with (\ref{pro3.12-a}) implies that
 $$
 \lim_{n\to\infty}\hat{E}[I_{\{|X_{n}-C|\geq \varepsilon\}}]=0,
 $$
i.e. $X_n\stackrel{V}{\longrightarrow} C.$\hfill\fbox

\subsection{Remarks}

\begin{rem}
An anonymous referee pointed out that quasi sure convergence does not imply convergence in capacity in $G$-expectation setting. Thus our assumption that  $\hat{E}$ has the monotone continuity property is strong in some sense. We would explore the corresponding problems under weaker conditions in future.
\end{rem}

\begin{rem} As to Slutsky's theorem under sublinear expectations, refer to \cite[Lemma 4.2]{Zhang15}.
\end{rem}

\begin{rem}
A natural question arise: can we give a Skorokhod type theorem under sublinear expectation? We tried, but we have not found the right way yet.
\end{rem}

\bigskip

{ \noindent {\bf\large Acknowledgments} \quad The authors  thank the anonymous referee for providing helpful comments to improve and clarify the manuscript.  This work was supported by National Natural Science Foundation of China (Grant No.
11371191).}

\end{document}